\input amstex\documentstyle{amsppt}  
\pagewidth{12.5cm}\pageheight{19cm}\magnification\magstep1
\topmatter
\title On the totally positive grassmannian\endtitle
\author G. Lusztig\endauthor
\address{Department of Mathematics, M.I.T., Cambridge, MA 02139}\endaddress
\thanks{Supported by NSF grant DMS-1566618.}\endthanks
\endtopmatter   
\document

\define\we{\wedge}

\define\lb{\linebreak}

\define\op{\oplus}
   
\define\part{\partial}
\define\emp{\emptyset}

\define\n{\notin}

\define\m{\mapsto}
\define\do{\dots}

\define\sub{\subset}

\define\nl{\newline}

\define\e{\epsilon}

\define\vp{\varpi}

\redefine\L{\Lambda}

\define\RR{\bold R}

\define\ZZ{\bold Z}

\define\cb{\Cal B}

\subhead 0.1\endsubhead
Let $V$ be an $\RR$-vector space of dimension $N\ge2$ with a fixed basis \lb $e_1,e_2,\do,e_N$.
Let $\cb$ be the manifold whose points are the complete flags
$F_1\sub F_2\sub\do\sub F_{N-1}$ in $V$. Here $F_i$ is a subspace of dimension $i$ of $V$.
In \cite{L94} I defined the totally positive part $\cb_{>0}$ of $\cb$, a certain open subset of $\cb$
 homeomorphic to $\RR_{>0}^{N(N-1)/2}$. We fix $k\in[1,N-1]$. 
Let $Gr^k$ be the manifold whose points are the subspaces of $V$ of dimension $k$. 
In \cite{L98} I defined the totally positive part $Gr^k_{>0}$ of $Gr^k$ as the image of $\cb_{>0}$ under the map
$\cb@>>>Gr^k$ which takes $F_1\sub F_2\sub\do\sub F_{N-1}$ to $F_k$. This is an open subset of $Gr^k$
homeomorphic to $\RR_{>0}^{k(N-k)}$. 
In \cite{L98, 3.4} it was shown that $Gr^k_{>0}$ can be described in terms of inequalities involving elements
in the canonical basis of the irreducible representation of $SL(V)$ corresponding to a multiple $c\vp_k$
of the $k$-th fundamental weight $\vp_k$ where $c$ is any integer $\ge N-1$. (The result in \cite{L98, 3.4} 
applies to any real partial flag manifold.) In a note added in the proof of \cite{L98} I stated 
(quoting Rietsch \cite{Ri}) that $c$ can be taken to be any number $\ge1$ (including $1$);
a similar statement was made for any partial flag manifold. However the proof in \cite{Ri} contained an error, 
see Geiss, Leclerc, Schr\"oer \cite{GLS}. (I thank B. Leclerc for providing this reference.)
In 2009, Rietsch (unpublished, but mentioned in \cite{KLS}) has shown that $c$ 
above can indeed be taken to be $1$.  Proofs of Rietsch's result appeared in \cite{TL} and \cite{La}.
But for a general partial flag 
manifold it is not known to what extent the result in \cite{L98, 3.4} can be improved. In this paper
we present a method which could possibly yield an improvement of the result of \cite{L98, 3.4}. In the case
of the Grassmannian this recovers the result of Rietsch, see Theorem 0.5; but one can hope that our method
applies also in other cases. This method is based on the observation of \cite{L98, \S2} that the positive part 
of a partial flag manifold is a single connected component of an explicit algebraic open subset of that partial 
flag manifold. For a further study of $Gr^k_{>0}$, see \cite{Po}. I thank P. Galashin and L. Williams for 
comments on an earlier version of this paper.
 
\subhead 0.2\endsubhead
{\it Notation.} For two integers $a\le b$ we set $[a,b]=\{z\in\ZZ;a\le z\le b\}$. 
For a finite set $I$ let $|I|$ be the cardinal of $I$.
For any $I\in[1,N]$ let $V_I$ be the subspace of $V$ with basis $\{e_i;i\in I\}$.
For any $I\sub[1,N]$ with $|I|=k$, we set $e_I=e_{i_1}\we e_{i_2}\we\do\we e_{i_k}\in\L^kV$
where $I$ consists of the numbers $i_1<i_2<\do<i_k$ and $\L^kV$ is the $k$-th exterior power of $V$.
Let $\L^kV_{>0}$ (resp. $\L^kV_{\ge0}$)
be the set of nonzero vectors in $\L^kV$ whose coordinates with respect to the basis
$\{e_I;|I|=k\}$ are all in $\RR_{>0}$ (resp. $\RR_{\ge0}$). 
Let $P\L^kV_{>0}$ (resp. $P\L^kV_{\ge0}$) be the set of lines in $\L^kV$ which are spanned by vectors in 
$\L^kV_{>0}$ (resp. $\L^kV_{\ge0}$). Define a linear map $A:V@>>>V$ by 
$A(e_1)=e_2$, $A(e_i)=e_{i+1}+e_{i-1}$ if $1<i<N$, $A(e_N)=e_{N-1}$. For any $r\in\RR_{\ge0}$ let 
$g_r=\exp(rA)\in GL(V)$. If $r>0$, the matrix of $g_r$ is totally positive.

\subhead 0.3\endsubhead
As in \cite{L98}, let $Gr^k_{\ge0}$ be the closure of $Gr^k_{>0}$ in $Gr^k$. 
We define $Gr'{}^k$ to be the set of all $E\in Gr^k$ such that $E\cap V_I=0$ for 
any $I\sub[1,N]$ such that $|I|=N-k$ (an open subset of $Gr^k$).
Let $Gr'{}^k_{>0}$ (resp. $Gr'{}^k_{\ge0}$) be the set of all $E\in Gr^k$ such that the line $\L^kE$ in 
$\L^kV$ is in $P\L^k(V)_{>0}$ (resp. $P\L^k(V)_{\ge0}$). According to \cite{L98, 3.2}, we have

(a) $Gr^k_{>0}\sub Gr'{}^k_{>0}$, $Gr^k_{\ge0}\sub Gr'{}^k_{\ge0}$.
\nl
We show:

(b) $Gr'{}^k_{>0}\sub Gr'{}^k$.
\nl
Assume that $E\in Gr'{}^k_{>0}$, $E\n Gr'{}^k$. We can find $I\sub[1,N]$ such that $|I|=N-k$, 
$E\cap V_I\ne0$. Let $e'_1\in E\cap V_I-\{0\}$.
We can find a basis $e'_1,e'_2,\do,e'_k$ of $E$ containing $e'_1$. Since $e'_1\in V_I$,
$\e:=e'_1\we e'_2\we\do\we e'_k$ is a linear combination of elements of the form $e_{I'}$
with $I'\sub[1,N],|I'|=k$, $I'\cap I\ne\emp$. In particular, $e_{[1,N]-I}$
appears with coefficient $0$ when $\e$ is expressed as a linear combination of $e_{I''},|I''|=k$. 
Thus, for any $c\in\RR-\{0\}$ we have $c\e\n\L^kV_{>0}$ so that $E\n Gr'{}^k$. This proves (b).

\subhead 0.4\endsubhead
Let $E\in Gr^k$. We say that $E$ is {\it generic} if

(i) $E\cap V_{[1,N-k]}=0$, 

(ii) $E\cap V_{[k+1,N]}=0$,

and if, setting

$E_i=E\cap V_{[1,N-k+i]}$ if $i\in[1,k-1]$, $E_k=E$, $E_i=E\op V_{[1,i-k]}$ if $i\in[k+1,N-1]$,

$E'_i=E\cap V_{[k-i+1,N]}$ if $i\in[1,k-1]$, $E'_k=E$, $E'_i=E\op V_{N-i+k+1,N}$ if 
$i\in[k+1,N-1]$,

so that $E_1\sub E_2\sub\do\sub E_{N-1}$, $\dim E_i=i$, $E'_1\sub E'_2\sub\do\sub E'_{N-1}$, 
$\dim E'_i=i$, we have:

(iii) $E'_i\cap E_{k-i}=0$ if $i\in[1,k-1]$,

(iv)  $E'_{k+i}\cap E_{N-i}=E$ if $i\in[1,N-k-1]$.
\nl
Let $Gr^{k*}$ be the set of all $E\in Gr_k$ which are generic. (An open subset of $Gr_k$.) 
According to \cite{L98}:

(a) {\it $Gr^k_{>0}$ is a connected component of $Gr^{k*}$.}
\nl
We show:

(b) $Gr'{}^k\sub Gr^{k*}$.
\nl
Let $E\in Gr'{}^k$. Then $E$ clearly satisfies conditions (i),(ii). For $i\in[1,k-1]$ we have
$$\align&E'_i\cap E_{k-i}=(E\cap V_{[k-i+1,N]})\cap(E\cap V_{[1,N-i]})\\&
=E\cap V_{[k-i+1,N]\cap[1,N-i]}=E\cap V_{[k-i+1,N-i]}=0\endalign$$
since $|[k-i+1,N-i]|=N-k$. Thus, (iii) holds. For $i\in[1,N-k-1]$ and for
$$x\in E'_{k+i}\cap E_{N-i}=(E\op V_{[N-i+1,N]})\cap(E\op V_{[1,N-k-i]})$$
we have $x=a+b=c+d$ with 
$$a\in E,b\in V_{[N-i+1,N]},c\in E,d\in V_{[1,N-k-i]}.$$
We have 
$$b-c\in V_{[N-i+1,N]}+V_{[1,N-k-i]}=V_{[N-i+1,N]\cup[1,N-k-i]}.$$
Also, $b-c\in E$ and $E\cap V_{[N-i+1,N]\cup [1,N-k-i]}=0$ since 
$|[N-i+1,N]\cup[1,N-k-i]|=N-k$. Thus $b=c$ so that $b=c=0$ since 
$V_{[N-i+1,N]}\cap V_{[1,N-k-i]}=0$. We see that $x=a\in E$. Thus, 
$E'_{k+i}\cap E_{N-i}\sub E$. The reverse inclusion is obvious. We see that (iv) holds. This proves (b).

\proclaim{Theorem 0.5 (Rietsch)} We have $Gr'{}^k_{>0}=Gr^k_{>0}$.
\endproclaim
The proof is similar to that of \cite{L94, 8.17}.
The inclusion $Gr^k_{>0}\sub Gr'{}^k_{>0}$ follows from 0.3(a). We show the reverse inclusion.
Let $E\in Gr'{}^k_{>0}$. Since $g_r$ is totally positive for $r>0$ and is the identity for $r=0$,
for any $r\in\RR_{\ge0}$ we have $g_r(E)\in Gr'{}^k_{>0}$. Applying \cite{L94, 5.2} to 
$g=g_1:\L^kV@>>>\L^kV$
and noting that $g_1^n=g_n$ for all integers $n\ge1$, we see that the sequence $g_n\L^kE$ ($n=1,2,\do$)
converges in $P\L^kV_{>0}$ to the Perron line $L_{g_1}\in P\L^kV_{>0}$ (notation of \cite{L94, 5.2(a)}).
From \cite{L94, 8.9(a)} we see that there exists $E_1\in Gr^k_{>0}$ such that $g_1(E_1)=E_1$. By 0.3(a) we 
have $E_1\in Gr'{}^k_{>0}$. 
Since $L_{g_1}$ is the unique $g_1$-stable line in $P\L^kV_{>0}$, we must have $L_{g_1}=\L^kE_1$.
Since $g_n\L^kE$ converges to $\L^kE_1$ as $n\to\infty$, it follows that $g_nE$
converges to $E_1$ as $n\to\infty$ (in $Gr^k$). Since 
$E_1\in Gr^k_{>0}$ and $Gr^k_{>0}$ is open in $Gr^k$, it follows that  
$g_{n_0}E\in Gr^k_{>0}$ for some integer $n_0\ge1$.
Since the map $r\m g_nE$ from $\RR_{\ge0}$ to $Gr'{}^k_{>0}$ is continuous, its image is contained in a
single connected component of $Gr^{k*}$ (recall that $Gr'{}^k_{>0}\sub Gr'{}^k\sub Gr^{k*}$, 
see 0.3(b),0.4(b)). In particular the image of
$0$ and that of $n_0$ (namely $E$ and $g_{n_0}E$) belong to the  
same connected component of $Gr^{k*}$. Since $g_{n_0}E\in Gr^k_{>0}$ and
$Gr^k_{>0}$ is a connected component of $Gr^{k*}$ (see 0.4(a)), it follows that $E\in Gr^k_{>0}$.
The theorem is proved.

\subhead 0.6\endsubhead
The following result is a consequence of the theorem in 0.5.

(a) $Gr'{}^k_{\ge0}=Gr^k_{\ge0}$.
\nl
The inclusion $Gr^k_{\ge0}\sub Gr'{}^k_{\ge0}$ follows from 0.3(a). We show the reverse inclusion.
Let $E\in Gr'{}^k_{\ge0}$.
Since $g_r$ is totally positive for $r>0$, for such $r$ we have  $g_r\L^kE\in P\L^kV_{>0}$ that is,
$\L^k(g_rE)\in P\L^kV_{>0}$. Using 0.5 we see that $g_rE\in Gr^k_{>0}$ for $r>0$. Taking the limit as
$r\to0$ we see that $E$ is in the closure of $Gr^k_{>0}$, that is, $E\in Gr^k_{\ge0}$. This proves (a).


\widestnumber\key{GLS}
\Refs
\ref\key{GLS}\by C. Geiss, B. Leclerc, 
J. Schr\"oer\paper Preprojective algebras and cluster algebras\inbook 
Trends in representation theory of algebras and related topics\bookinfo EMS Ser. Congr. Rep.,
Eur. Math. Soc.,Z\"urich\yr2008\pages253-283\endref
\ref\key{KLS}\by A. Knutson, T. Lam, D. Speyer\paper Positroid varieties I. Juggling and geometry\lb \jour
arxiv:0903.3694\endref
\ref\key{La}\by T. Lam\paper Totally nonnegative Grassmannian and Grassmann polytopes
\jour arxiv:1506.00603\endref
\ref\key{L94}\by G. Lusztig\paper Total positivity in reductive groups\inbook Lie theory and geometry\bookinfo
 Progr. in Math. 123\publ Birkh\"auser Boston \yr1994\pages 531-568\endref
\ref\key{L98}\by G. Lusztig\paper Total positivity in partial flag manifolds\jour Represent. Th.\vol2\yr1998
\pages 70-78\endref
\ref\key{Po}\by A. Postnikov\paper Total positivity, Grassmannians and networks\jour arxiv:math.0609764\endref
\ref\key{Ri}\by K. Rietsch\paper Total positivity and real flag varieties\jour Ph.D. Thesis, MIT\yr1998\endref
\ref\key{TL}\by K. Talaska, L. Williams\paper Network parametrizations for the grassmannian\jour arxiv:1210.5433
\endref
\endRefs
\enddocument